\documentclass[12pt,a4paper]{article}
\usepackage{indentfirst}
\usepackage{amsmath,amssymb,amsfonts,amsthm,graphics}
\usepackage{makeidx}
\usepackage{amsmath,amssymb,amsfonts,amsthm,graphics,graphicx}
\usepackage{makeidx}
\usepackage{color}
\usepackage{ifpdf}
\usepackage{subfigure}
\usepackage{mathrsfs}

\newtheorem{thm}{Theorem}

\newtheorem{prob}{Problem}
\newtheorem{lem}{Lemma}
\newtheorem{pro}{Proposition}

\theoremstyle{definition}

\def\-{\mbox{--}}

\newtheorem{remark}{Remark}

\newtheorem{conj}{Conjecture}

\newtheorem{ex}{Example}

\def\pf{\noindent {\it Proof.} }
\begin{document}

\title{\Large\bf Total monochromatic connection of graphs\footnote{Supported by NSFC No.11371205 and 11531011, and PCSIRT.} }
\author{\small Hui Jiang, Xueliang Li, Yingying Zhang\\
\small Center for Combinatorics and LPMC\\
\small Nankai University, Tianjin 300071, China\\
\small E-mail: jhuink@163.com; lxl@nankai.edu.cn;\\
\small zyydlwyx@163.com}
\date{}
\maketitle
\begin{abstract}

A graph is said to be {\it total-colored} if all the edges and the vertices of the graph are colored. A path in a total-colored graph is a {\it total monochromatic path} if all the edges and internal vertices on the path have the same color. A total-coloring of a graph is a {\it total monochromatically-connecting coloring} ({\it TMC-coloring}, for short) if any two vertices of the graph are connected by a total monochromatic path of the graph. For a connected graph $G$, the {\it total monochromatic connection number}, denoted by $tmc(G)$, is defined as the maximum number of colors used in a TMC-coloring of $G$. These concepts are inspired by the concepts of monochromatic connection number $mc(G)$, monochromatic vertex connection number $mvc(G)$ and total rainbow connection number $trc(G)$ of a connected graph $G$. Let $l(T)$ denote the number of leaves of a tree $T$, and let $l(G)=\max\{ l(T) | $ $T$ is a spanning tree of $G$ $\}$ for a connected graph $G$. In this paper, we show that there are many graphs $G$ such that $tmc(G)=m-n+2+l(G)$, and moreover, we prove that for almost all graphs $G$, $tmc(G)=m-n+2+l(G)$ holds. Furthermore, we compare $tmc(G)$ with $mvc(G)$ and $mc(G)$, respectively, and obtain that there exist graphs $G$ such that $tmc(G)$ is not less than $mvc(G)$ and vice versa, and that $tmc(G)=mc(G)+l(G)$ holds for almost all graphs. Finally, we prove that $tmc(G)\leq mc(G)+mvc(G)$, and the equality holds if and only if $G$ is a complete graph.


\end{abstract}

\section{Introduction}

In this paper, all graphs are simple, finite and undirected. We refer to the book \cite{B} for undefined notation and terminology in graph theory. Throughout this paper, let $n$ and $m$ denote the order (number of vertices) and size (number of edges) of a graph, respectively. Moreover, a vertex of a connected graph is called a {\it leaf} if its degree is one; otherwise, it is an {\it internal vertex}. Let $l(T)$ and $q(T)$ denote the number of leaves and the number of internal vertices of a tree $T$, respectively, and let $l(G)=\max\{ l(T) | $ $T$ is a spanning tree of $G$ $\}$ and $q(G)=\min\{ q(T) | $ $T$ is a spanning tree of $G$ $\}$ for a connected graph $G$. Note that the sum of $l(G)$ and $q(G)$ is $n$ for any connected graph $G$ of order $n$. A path in an edge-colored graph is a {\it monochromatic path} if all the edges on the path have the same color. An edge-coloring of a connected graph is a {\it monochromatically-connecting coloring} ({\it MC-coloring}, for short) if any two vertices of the graph are connected by a monochromatic path of the graph. For a connected graph $G$, the {\it monochromatic connection number} of $G$, denoted by $mc(G)$, is defined as the maximum number of colors used in an MC-coloring of $G$. An {\it extremal MC-coloring} is an MC-coloring that uses $mc(G)$ colors. Note that $mc(G)=m$ if and only if $G$ is a complete graph. The concept of $mc(G)$ was first introduced by Caro and Yuster \cite{CY} and has been well-studied recently. We refer the reader to \cite{CLW,GLQZ} for more details.

As a natural counterpart of the concept of monochromatic connection, Cai et al. \cite{CLW1} introduced the concept of monochromatic vertex connection. A path in a vertex-colored graph is a {\it vertex-monochromatic path} if its internal vertices have the same color. A vertex-coloring of a graph is a {\it monochromatically-vertex-connecting coloring} ({\it MVC-coloring}, for short) if any two vertices of the graph are connected by a vertex-monochromatic path of the graph. For a connected graph $G$, the {\it monochromatic vertex connection number}, denoted by $mvc(G)$, is defined as the maximum number of colors used in an MVC-coloring of $G$. An {\it extremal MVC-coloring} is an MVC-coloring that uses $mvc(G)$ colors. Note that $mvc(G)=n$ if and only if $diam(G)\leq 2$.

Actually, the concepts of monochromatic connection number $mc(G)$ and monochromatic vertex connection number $mvc(G)$ are natural opposite concepts of rainbow connection number $rc(G)$ and rainbow vertex connection number $rvc(G)$. For details about them we refer to a book \cite{LiS} and a survey paper \cite{LiSS}. Moreover, the concept of total rainbow connection number $trc(G)$ in \cite{LMS1} was motivated by the rainbow connection number $rc(G)$ and rainbow vertex connection number $rvc(G)$. Thus, here we introduce the concept of total monochromatic connection of graphs. A graph is said to be {\it total-colored} if all the edges and the vertices of the graph are colored. A path in a total-colored graph is a {\it total monochromatic path} if all the edges and internal vertices on the path have the same color. A total-coloring of a graph is a {\it total monochromatically-connecting coloring} ({\it TMC-coloring}, for short) if any two vertices of the graph are connected by a total monochromatic path of the graph. For a connected graph $G$, the {\it total monochromatic connection number}, denoted by $tmc(G)$, is defined as the maximum number of colors used in a TMC-coloring of $G$. An {\it extremal TMC-coloring} is a TMC-coloring that uses $tmc(G)$ colors. It is easy to check that $tmc(G)=m+n$ if and only if $G$ is a complete graph.

The rest of this paper is organized as follows: In Section $2$, we prove that $tmc(G)\geq m-n+2+l(G)$ for any connected graph and determine the value of $tmc(G)$ for some special graphs. In Section $3$, we prove that there are many graphs with $tmc(G)=m-n+2+l(G)$ which are restricted by other graph parameters such as the maximum degree, the diameter and so on, and moreover, we show that for almost all graphs $G$, $tmc(G)=m-n+2+l(G)$ holds. In Section $4$, we compare $tmc(G)$ with $mvc(G)$ and $mc(G)$, respectively, and obtain that there exist graphs $G$ such that $tmc(G)$ is not less than $mvc(G)$ and vice versa, and that $tmc(G)=mc(G)+l(G)$ for almost all graphs. Moreover, we prove that $tmc(G)\leq mc(G)+mvc(G)$, and the equality holds if and only if $G$ is a complete graph.

\section{Preliminary results}

In this section, we show that $tmc(G)\geq m-n+2+l(G)$ and present some preliminary results on the total monochromatic connection number. Moreover, we determine the value of $tmc(G)$ when $G$ is a tree, a wheel, and a complete multipartite graph. The following fact is easily seen.

\begin{pro}\label{pro1} If $G$ is a connected graph and $H$ is a connected spanning subgraph of $G$, then $tmc(G)\geq e(G)-e(H)+tmc(H)$.
\end{pro}

Since for any two vertices of a tree, there exists only one path connecting them, we have the following result.

\begin{pro}\label{pro2} If $T$ is a tree, then $tmc(T)=l(T)+1$.
\end{pro}

The consequence below is immediate from Propositions \ref{pro1} and \ref{pro2}.

\begin{thm}\label{thm1} For a connected graph $G$, $tmc(G)\geq m-n+2+l(G)$.
\end{thm}

Let $G$ be a connected graph and $f$ be an extremal TMC-coloring of $G$ that uses a given color $c$. Note that the subgraph $H$ formed by the edges and vertices colored $c$ is connected, or we will give a fresh color to all the edges and vertices colored $c$ in some of these components while still maintaining a TMC-coloring. Moreover, the color of each internal vertex of $H$ is $c$. Otherwise, let $u_1,\ldots,u_t$ be the internal vertices of $H$ such that each of them is not colored $c$. We obtain the subgraph $H_0$ by deleting the vertices $\{u_1,\ldots,u_t\}$. If $H_0$ is connected, it is possible to choose an edge incident with $u_1$ and assign it with a fresh color while still maintaining a TMC-coloring. If not, we can give a fresh color to all the edges and vertices colored $c$ in some of these components while still maintaining a TMC-coloring. Furthermore, $H$ does not contain any cycle; otherwise, a fresh color can be assigned to any edge of the cycle while still maintaining a TMC-coloring. Thus, $H$ is a tree where the color of each internal vertex is $c$. Now we define the {\it color tree} as the tree formed by the edges and vertices colored $c$, denoted by $T_c$. If $T_c$ has at least two edges, the color $c$ is called {\it nontrivial}. Otherwise, $c$ is {\it trivial}. We call an extremal TMC-coloring {\it simple} if for any two nontrivial colors $c$ and $d$, the corresponding trees $T_c$ and $T_d$ intersect in at most one vertex. If $f$ is simple, then the leaves of $T_c$ must have distinct colors different from color $c$. Otherwise, we can give a fresh color to such a leaf while still maintaining a TMC-coloring. Moreover, a nontrivial color tree of $f$ with $m'$ edges and $q'$ internal vertices is said to {\it waste} $m'-1+q'$ colors. For the rest of this paper we will use these facts without further mentioning them.

The lemma below shows that one can always find a simple extremal TMC-coloring for a connected graph.

\begin{lem}\label{lem1} Every connected graph $G$ has a simple extremal TMC-coloring.
\end{lem}

\pf We are given an extremal TMC-coloring $f$ of $G$ with the most number of trivial colors, and then we prove that this coloring must be simple. Suppose that there exist two nontrivial colors $c$ and $d$ such that $T_c$ and $T_d$ contain $k$ common vertices denoted by $u_1,u_2,\ldots,u_k$, where $k\geq 2$. Now we divide our discussion into two cases.

\textbf{Case 1.} For $1\leq i\leq k$, $u_i$ is an internal vertex of $T_c$ or $T_d$.

For $1\leq i\leq k$, if $u_i$ is an internal vertex of $T_c$, $u_i$ must be a leaf of $T_d$ and then set $e_i=u_iw_i$ where $w_i$ is the neighbor of $u_i$ in $T_d$; otherwise, $u_i$ must be a leaf of $T_c$ and then put $e_i=u_iv_i$ where $v_i$ is the neighbor of $u_i$ in $T_c$. Let $H$  denote the subgraph consisting of the edges and vertices of $T_c\cup T_d$. Clearly, $H$ is connected. We obtain a spanning tree $H_0$ of $H$ by deleting the edges $\{e_2,e_3,\ldots,e_k\}$. Now we change the total-coloring of $H$ while still maintaining the colors of the leaves in $H_0$ unchanged. Assign the edges and internal vertices of $H_0$ with color $c$ and the remaining edges $\{e_2,e_3,\ldots,e_k\}$ with distinct new colors. Obviously, the new total-coloring is also a TMC-coloring and uses $k-2$ more colors than our original one. So, it either uses more colors or uses the same number of colors but more trivial colors, contradicting the assumption on $f$.

\textbf{Case 2.} There exists a vertex among $u_1,\ldots,u_k$, say $u_1$, which is a leaf of both $T_c$ and $T_d$.

Let $v_1$ and $w_1$ be the neighbors of $u_1$ in $T_c$ and $T_d$, respectively. There must be another color tree $T_e$ (including a single edge) connecting $v_1$ and $w_1$. For $1\leq i\leq k$, if $u_i$ is a leaf of $T_c$, then set $e_i=u_iv_i$ where $v_i$ is the neighbor of $u_i$ in $T_c$; otherwise, $u_i$ must be a leaf of $T_d$ and then put $e_i=u_iw_i$ where $w_i$ is the neighbor of $u_i$ in $T_d$. Let $H_1$ denote the subgraph consisting of the edges and vertices of $T_c\cup T_d$. We obtain a spanning subgraph $H_2$ of $H_1$ by deleting the edges $\{e_1,e_2,\ldots,e_k\}$. If $T_e$ and $H_2$ do not have common leaves, let $E_0=\{e_1,e_2,\ldots,e_k\}$. Otherwise, let $u'_1,\ldots,u'_t$ denote the common leaves of $T_e$ and $H_2$. Set $e'_i=u'_iv'_i$ where $v'_i$ is the neighbor of $u'_i$ in $T_e$ for $1\leq i\leq t$. And then let $E_0=\{e_1,\ldots,e_k,e'_1,\ldots,e'_t\}$.
Let $H$ denote the subgraph  consisting of the edges and vertices of $T_c\cup T_d\cup T_e$. Clearly, $H$ is connected. We obtain a spanning connected subgraph $H_0$ of $H$ by deleting the edges of $E_0$. Now we change the total-coloring of $H$ while still maintaining the colors of the leaves in $H_0$ unchanged. Assign the edges and internal vertices of $H_0$ with color $c$ and the remaining edges of $H$ (i.e., the edges of $E_0$) with distinct new colors. Note that if $v$ is a common leaf of either $T_c$ and $T_d$ or $T_e$ and $H_2$, it is also a leaf of $H_0$. Obviously, the new total-coloring is also a TMC-coloring and uses at least $k+t-2$ more colors than our original one. So, it either uses more colors or uses the same number of colors but more trivial colors, contradicting the assumption on $f$.
 \qed

Now we use the above results to compute the total monochromatic connection numbers of wheel graphs and complete multipartite graphs.

\begin{ex}\label{ex1} Let $G$ be a wheel $W_{n-1}$ of order $n\geq 5$. Then $tmc(G)=m-n+2+l(G)$.
\end{ex}

\pf We are given a simple extremal TMC-coloring $f$ of $G$. Note that $m-n+2+l(G)=m+1$ and $tmc(G)\geq m+1$ by Theorem \ref{thm1}. Suppose that $f$ consists of $k$ nontrivial color trees, denoted by $T_1,\ldots,T_k$. In fact, we can always find two vertices with degree at least 4 if $k\geq 3$, a contradiction. Likewise, if $k=2$, $G$ must be $W_4$ and $tmc(W_4)=m+1$. Thus, assume that $k=1$ and $T_1$ is not spanning (Otherwise, $tmc(G)=m-n+2+l(G)$). Note that for every vertex $v\notin T_1$, there exist the total monochromatic paths connecting $v$ and the $|T_1|$ vertices of $T_1$. As $f$ is simple, these paths are internally vertex-disjoint. Hence, $deg(v)\geq |T_1|$. If $|T_1|\geq 4$, the $n-1$ vertices with degree 3 of $G$ must be in $T_1$ and then $T_1$ is a path. Thus, $tmc(G)= m+n-(n-3)-(n-3)=m+6-n\leq m+1$. If $|T_1|=3$, then $G$ must be $W_3$ while $n\geq 5$. Therefore, the proof is complete.
\qed

\begin{ex}\label{ex2}Let $G= K_{n_1,\ldots,n_r}$ be a complete multipartite graph with $n_1 \geq \ldots \geq n_t\geq 2$ and $n_{t+1}=\ldots=n_r=1$. Then $tmc(G)=m+r-t$.
\end{ex}

\pf The case that $r=2$ is a special case of Theorem \ref{thm3} whose proof is given in Section 3, so assume that $r\geq 3$. Let $f$ be a simple extremal TMC-coloring of $G$ with maximum trivial colors. Suppose that $f$ consists of $k$ nontrivial color trees, denoted by $T_1,\ldots,T_k$, where $t_i=|V(T_i)|$ and $q_i=q(T_i)$ for $1\leq i\leq k$. Now we divide our discussion into two cases.

\textbf{Case 1.} $t=r$.

In this case, every vertex appears in at least one of the nontrivial color trees. Note that $m-n+2+l(G)=m$ and $tmc(G)\geq m$ by Theorem \ref{thm1}. If $\sum_{i=1}^k(t_i-1)\geq n$, then we have that $tmc(G)\leq m+n-n-\sum_{i=1}^kq_i+k =m-\sum_{i=1}^kq_i+k\leq m$. Thus, $tmc(G)=m$. Suppose that $\sum_{i=1}^k(t_i-1)\leq n-1$. Now consider the subgraph $G'$ consisting of the union of the $T_i$ and let $C_1,\ldots,C_s$ denote its components.

Now we may assume that there exists a component, say $C_1$, such that each nontrivial color tree in $C_1$ is a star. Let $S$ be a star of $C_1$ with center $u$ and leaves $u_1,\ldots,u_p$, where $u_1,\ldots,u_{p'}$ are in the same vertex class, say $V_1$. Suppose that $p'\geq 2$. Indeed, if $p'=1$, we can give a new color to the edge $uu_1$ while still maintaining a TMC-coloring. We claim that $C_1$ contains a cycle. If $p'<|V_1|$, there exists a vertex $u_{p+1}$ of $V_1$ not adjacent to $u$ in $S$. Then $u_1$ and $u_{p+1}$ must be in a same nontrivial color tree and the same happens for $u_{p'}$ and $u_{p+1}$. These nontrivial color trees containing $u_1$, $u_{p'}$ and $u_{p+1}$ must form a cycle. If $p'=|V_1|$, we have that the vertices of the vertex class containing $u$ must be in a same nontrivial color tree, or we will get a cycle in a similar way. By that analogy, we obtain a cycle formed by some centers of the nontrivial color trees in $C_1$. Now we change the total-coloring of $C_1$. We obtain a spanning tree $T'$ of $C_1$ by connecting $u_1$ to the vertices in the same class with $u$ and $u$ to the other vertices of $C_1$. We color the edges and internal vertices of $T'$ with the same color and all other edges and vertices with distinct new colors. Clearly, this new total-coloring is also a TMC-coloring. However, it either uses more colors or uses the same number of colors but more trivial colors, contradicting the assumption on $f$.

Thus, suppose that there exists a nontrivial color tree of $C_i$, say $T_{i1}$, having two adjacent internal vertices $u_i$ and $v_i$ for $1\leq i\leq s$. We obtain a spanning tree $T$ by connecting $v_1$ to each vertex in the same class with $u_1$ of $G$ and $u_1$ to the other vertices of $G$. Now we give a new total-coloring $f'$ of $G$. Color the edges and internal vertices of $T$ with the same color and all other edges and vertices of $G$ with distinct new colors. Obviously, $f'$ is still a TMC-coloring. If $s\geq 2$, then it either uses more colors or uses the same number of colors but more trivial colors than $f$, a contradiction. Thus, $s=1$. Moreover, we can check that $f'$ is a simple extremal TMC-coloring with maximum trivial colors. Therefore, $tmc(G)=m$.

\textbf{Case 2.} $t<r$.

We obtain a star $S$ by connecting a vertex of $\cup_{i=t+1}^{r}V_i$ to each vertex of $\cup_{i=1}^t V_i$. Color the edges and the center vertex of $S$ with the same color and all other edges and vertices of $G$ with distinct new colors. Clearly, this new total-coloring is still a TMC-coloring, denoted by $f'$. Thus, $tmc(G)\geq m+r-t$. If $\sum_{i=1}^k(t_i-1)\geq n-r+t$, then we have that
$tmc(G)\leq m+n-(n-r+t)-\sum_{i=1}^kq_i+k=m+r-t-
\sum_{i=1}^k q_i+k\leq m+r-t$. Hence, $tmc(G)= m+r-t$. Suppose that $\sum_{i=1}^k(t_i-1)\leq n-r+t-1$. Next consider the subgraph $G'$ consisting of the union of the $T_i$'s and suppose that it has $s$ components, say $C_1,\ldots,C_s$. Note that $|V(G')|\geq n-r+t$ since any two vertices of the same class must be covered in a nontrivial color tree. The case that $|V(G')|= n-r+t$ can be verified by a similar discussion to Case 1. Thus, suppose that $|V(G')|>n-r+t$. It is obvious that $s\geq 2$. Moreover, there must exist a vertex $x$ of $\cup_{i=t+1}^r V_i$, which is contained in a component of $G'$, say $C_1$. For $2\leq j\leq s$, there does not exist a vertex of $\cup_{i=t+1}^r V_i$ in $C_j$. Otherwise, let $x$ be the center of $S$ and then $f'$ either uses more colors or uses the same number of colors but more trivial colors than $f$, a contradiction. By a similar discussion to Case 1, we can obtain that there exists a nontrivial color tree of $C_j$ having two adjacent internal vertices for $2\leq j\leq s$. We obtain a star $S_1$ by joining the vertices of $\cup_{i=2}^{s}C_i$ to one internal vertex of $C_1$. We give a new total-coloring of $G$ while still maintaining the total-coloring of $C_1$ unchanged. Assign the edges and the center vertex of $S_1$ with one color and the other edges and vertices of $G\backslash C_1$ with distinct new colors. This new total-coloring is still a TMC-coloring and it either uses more colors or uses the same number of colors but more trivial colors, contradicting the assumption on $f$. Therefore, we have finished the proof.
\qed

\section{Graphs with $tmc(G)=m-n+2+l(G)$}

In this section, we prove that there are many graphs $G$ for which $tmc(G)=m-n+2+l(G)$, even for almost all graphs.

\begin{lem}\label{lem2}\cite{CY} Let $G$ be a connected graph of order $n>3$. If $G$ satisfies any of the following properties, then $mc(G)=m-n+2$.

$(a)$ The complement $\overline{G}$ of $G$ is $4$-connected.

$(b)$ $G$ is $K_3$-free.

$(c)$ $\Delta(G)<n-\frac{2m-3(n-1)}{n-3}$. In particular, this holds if $\Delta(G)\leq(n+1)/2$, and this also holds if $\Delta(G)\leq n-2m/n$.

$(d)$ $diam(G)\geq 3$.

$(c)$ $G$ has a cut vertex.
\end{lem}

We can obtain that $tmc(G)\leq mc(G)+l(G)$ for a noncomplete graph, whose proof is contained in the proof of Theorem \ref{thm6} in Section 4. In addition with Theorem \ref{thm1} and Lemma \ref{lem2}, we have the following results.

\begin{thm}\label{thm2} Let $G$ be a connected graph of order $n>3$. If $G$ satisfies any of the following properties, then $tmc(G)=m-n+2+l(G)$.

$(a)$ The complement $\overline{G}$ of $G$ is $4$-connected.

$(b)$ $G$ is $K_3$-free.

$(c)$ $\Delta(G)<n-\frac{2m-3(n-1)}{n-3}$.

$(d)$ $diam(G)\geq 3$.

$(c)$ $G$ has a cut vertex.
\end{thm}

One cannot hope to strengthen Theorem \ref{thm2}(c) by improving the upper bound of $\Delta(G)$. In fact, let $G=K_{n-2,1,1}$. Then we have that $tmc(G)=m-n+3+l(G)$ while the maximum degree is $n-1=n-\frac{2m-3(n-1)}{n-3}$.

From Theorem \ref{thm2}(a), we can get a stronger result. For a property $P$ of graphs and a positive integer $n$, define $Prob(P,n)$ to be the ratio of the number of graphs with $n$ labeled vertices having $P$ over the total number of graphs with these vertices. If $Prob(P,n)$ approaches 1 as $n$ tends to infinity, then we say that {\it almost all} graphs have the property $P$. See \cite{BH} for example.

\begin{thm}\label{thm3} For almost all graphs $G$, we have that $tmc(G)=m-n+2+l(G)$.
\end{thm}

In order to prove Theorem \ref{thm3}, we need the following lemma.

\begin{lem}\cite{BH}\label{lem3}  For every nonnegative integer $k$, almost all graphs are $k$-connected.
\end{lem}

\noindent {\it Proof of Theorem 3:}  For any given nonnegative integer $n$, let $\mathscr{G}_n$ denote the set of all graphs of order $n$, and let $\mathscr{G}_n^4$ denote the set of all 4-connected graphs of order $n$. Moreover, let $\mathscr{B}_n$ denote the set of all graphs $G$ of order $n$ such that the complement $\overline{G}$ of $G$ is 4-connected. Note that for any two graphs $G$ and $H$, $G\cong H$ if and only if $\overline{G}\cong\overline{H}$. Then, it is easy to check that the map: $G\rightarrow\overline{G}$ is a bijection from $\mathscr{B}_n$ to $\mathscr{G}_n^4$. Therefore, we have
$$\frac{|\mathscr{B}_n|}{|\mathscr{G}_n|}=\frac{|\mathscr{G}_n^4|}{|\mathscr{G}_n|}.$$
By Lemma \ref{lem3}, it follows that almost all graphs are 4-connected. Then, we get that almost all graphs have 4-connected complements. Furthermore, since almost all graphs are connected, we have that $tmc(G)=m-n+2+l(G)$ by Theorem \ref{thm2}(a).
\qed

\begin{remark}\label{rem1} For the monochromatic connection number $mc(G)$, from Lemma \ref{lem2}(a) and Lemma \ref{lem3}, one can deduce, in a similar way, that for almost all graphs $G$, $mc(G)=m-n+2$ holds.
\end{remark}
\begin{remark}\label{rem2}
To use the parameter $l(G)$ in the above formulas looks good. However, from \cite[p.206]{GS} we know that it is NP-hard to find a spanning tree that has maximum number of leaves in a connected graph $G$.
\end{remark}

\section{Compare $tmc(G)$ with $mvc(G)$ and $mc(G)$ }

Let $G$ be a nontrivial connected graph. Firstly, we compare $tmc(G)$ with $mvc(G)$. The question we may ask is, can we bound one of $tmc(G)$ and $mvc(G)$ in terms of the other? The following two theorems give sufficient conditions for $tmc(G)>mvc(G)$.

\begin{thm}\label{thm4} Let $G$ be a connected graph with diameter $d$. If $m\geq 2n-d-2$, then $tmc(G)>mvc(G)$.
\end{thm}

\pf The case that $d=1$ is trivial, so assume that $d\geq 2$. We can check that if $l(G)=2$, then $tmc(G)>mvc(G)$. Thus, suppose that $l(G)\geq 3$. By Theorem \ref{thm1}, it follows that $tmc(G)\geq m-n+2+l(G)\geq 2n-d-2-n+2+3=n-d+3$. Moreover, we have that $mvc(G)\leq n-d+2$ by \cite[Proposition 2.3]{CLW1}. Therefore, $tmc(G)>mvc(G)$.
\qed

\begin{thm}\label{thm5} Let $G$ be a connected graph of diameter $2$ with maximum degree $\Delta$. If $\Delta\geq \frac{n+1}{2}$, then $tmc(G)>mvc(G)$.
\end{thm}

Before proving Theorem \ref{thm5}, we need the lemma below.

\begin{lem}\cite{BO}\label{lem4} Let $G$ be a connected graph of diameter $2$ with maximum degree $\Delta$. Then
\begin{eqnarray}m\geq
\begin{cases}
n+\Delta-2, &if\ \Delta=n-2\ or\ n-3\cr
2n-5, &if\ \Delta=n-4\cr
2n-4, &if\ \frac{2n-2}{3}\leq \Delta\leq n-5\cr
3n-\Delta-6, &if\ \frac{3n-3}{5}\leq \Delta< \frac{2n-2}{3}\cr
5n-4\Delta-10, &if\ \frac{5n-3}{9}\leq \Delta< \frac{3n-3}{5}\cr
4n-2\Delta-11, &if\ \frac{n+1}{2}\leq \Delta< \frac{5n-3}{9}
\end{cases}
\end{eqnarray}
\end{lem}

\noindent {\it Proof of Theorem 5:} The case that $n\leq 7$ can be easily verified. Suppose that $n\geq 8$. Since the diameter of $G$ is $2$, we have that $mvc(G)=n$. By Theorem \ref{thm1} and Lemma \ref{lem4}, $tmc(G)\geq m-n+2+l(G)>n$. Thus, $tmc(G)>mvc(G)$.
\qed

Actually, we have that $tmc(C_5)=4<mvc(C_5)=5$, where $m<2n-d-2$ and $\Delta < \frac{n+1}{2}$. This implies that the conditions of Theorems \ref{thm4} and \ref{thm5} cannot be improved. Moreover, if $G$ is a star, then $tmc(G)=mvc(G)=n$. Therefore, there exist graphs $G$ such that $tmc(G)$ is not less than $mvc(G)$ and vice versa. However, we cannot show whether there exist other graphs with $tmc(G)\leq mvc(G)$. Thus, we propose the following problem.

\begin{prob}\label{prob1} Dose there exist a graph of order $n\geq 6$ except a star such that $tmc(G)\leq mvc(G)$?
\end{prob}

Next we compare $tmc(G)$ with $mc(G)$. If $G$ satisfies one of the conditions in Theorem \ref{thm2}, then we have $mc(G)=m-n+2$ and so $tmc(G)=mc(G)+l(G)$. For a complete graph $G$, $tmc(G)>mc(G)+l(G)$. From \cite[Corollary 13]{CY}, if $G$ is a wheel $W_{n-1}$ of order $n\geq 5$, we have that $mc(G)=m-n+3$ and then $tmc(G)<mc(G)+l(G)$.
However, by Theorem \ref{thm3} and Remark \ref{rem1}, it follows that almost all graphs have that $tmc(G)=mc(G)+l(G)$ which implies that almost all graphs have that $tmc(G)>mc(G)$. Thus, we propose the following conjecture.

\begin{conj}\label{conj1} For a connected graph $G$, it always holds that $tmc(G)>mc(G)$.
\end{conj}

Finally, we compare $tmc(G)$ with $mc(G)+mvc(G)$.

\begin{thm}\label{thm6} Let $G$ be a connected graph. Then $tmc(G)\leq mc(G)+mvc(G)$, and the equality holds if and only if $G$ is a complete graph.
\end{thm}

In order to prove Theorem \ref{thm6}, we need the following lemma.

\begin{lem}\label{lem5} For a noncomplete connected graph $G$, let $f$ be a simple extremal TMC-coloring of $G$ and $T_1,\ldots,T_k$ denote all the nontrivial color trees of $f$, where $t_i=|V(T_i)|$ and $q_i=q(T_i)$ for $1\leq i\leq k$. Then, $\sum_{i=1}^{k}q_i\geq q(G)$.
\end{lem}

\pf For any $v\in G$, if $v\notin \cup_{i=1}^{k}T_i$, $v$ must be adjacent to an internal vertex $w_0$ of a nontrivial color tree and then set $E_v=\{vw|w\in N(v)\backslash\{w_0\}\}$. If $v$ is an internal vertex of a nontrivial color tree containing $v$, set $E_v=\emptyset$. Otherwise, $v$ is a leaf of any nontrivial color tree containing $v$. Let $T_1,\ldots,T_s$ denote the nontrivial color trees containing $v$ and $v_1,\ldots,v_s$ be the neighbors of $v$ in $T_1,\ldots,T_s$, respectively. Let $E_v=\{vv_2,\ldots,vv_s\}$. We obtain a spanning subgraph $G'$ by deleting the edges of $\bigcup_{v\in G}E_v$. Note that every vertex of $\{v:E_v=\emptyset\}$ is connected to each other. For any two vertices $u_1$ and $u_2$ of $\{v:E_v=\emptyset\}$, there exists a total monochromatic path $P$ of $G$ connecting them. For each vertex $u$ of $P$, we have $E_u=\emptyset$. Thus, $G'$ also contains $P$ from $u_1$ to $u_2$. Moreover, every vertex of $\{v:E_v\neq\emptyset\}$ is connected to a vertex of $\{v:E_v=\emptyset\}$. Hence, $G'$ is connected and each vertex of $\{v:E_v\neq\emptyset\}$ cannot be an internal vertex of $G'$. Then $\sum_{i=1}^{k}q_i\geq q(G')\geq q(G)$.
\qed

Now, we are ready to prove Theorem \ref{thm6}.

\noindent {\it Proof of Theorem 6:} If $G$ is a complete graph, we have that $tmc(G)=mc(G)+mvc(G)$. Thus, suppose that $G$ is not complete. We are given a simple extremal TMC-coloring $f$ of $G$. Suppose that $f$ consists of $k$ nontrivial color trees denoted by $T_1,\ldots,T_k$, where $t_i=|V(T_i)|$ and $q_i=q(T_i)$ for $1\leq i\leq k$. Then $tmc(G)=m+n-\sum_{i=1}^{k}(t_i-2)-\sum_{i=1}^{k}q_i$. Now we take a copy $G'$ of $G$. Then $G'$ contains the trees $T'_1,\ldots,T'_k$ corresponding to $T_1,\ldots,T_k$, respectively. Define an edge-coloring $f_e$ of $G'$ as follows: color the edges of $T_i$ with color $i$ for $1\leq i\leq k$ and the other edges of $G'$ with distinct new colors. Then $f_e$ is an MC-coloring of $G'$ with $m-\sum_{i=1}^{k}(t_i-2)$ colors. Thus, $mc(G)=mc(G')\geq m-\sum_{i=1}^{k}(t_i-2)=tmc(G)-n+\sum_{i=1}^{k}q_i$. By Lemma \ref{lem5}, we have that $\sum_{i=1}^{k}q_i\geq q(G)$. Then $tmc(G)\leq mc(G)+n-q(G)=mc(G)+l(G)$. Moreover, it is easy to obtain that $mvc(G)\geq l(G)+1$. Hence, $tmc(G)<mc(G)+mvc(G)$. Therefore, the proof is complete.
\qed

\begin{remark}\label{rem3} For the total rainbow connection number $trc(G)$, we cannot bound one of $trc(G)$ and $rc(G)+rvc(G)$ in terms of the other. For a connected graph $G$, $trc(G)=rc(G)+rvc(G)$ if $G$ is a complete graph or a star. Moreover, if $G$ is a complete bipartite graph $K_{m,n}$ with $m\geq 2$ and $n\geq 6^{m}$, then $trc(G)=7>rc(G)+rvc(G)=4+1$ \cite{LiSS,LMS,LMS1}. In \cite{LMS1}, for every $s\geq 1481$, there exists a graph $G$ with $trc(G)=rvc(G)=s$ which implies that $trc(G)<rc(G)+rvc(G)$. This is one thing that the total monochromatic connection differs from the total rainbow connection.
\end{remark}

\end{document}